\documentclass[10pt,twoside]{classecras-e}
\usepackage{amssymb,amsbsy,amsmath,amsfonts,amssymb,amscd}
\usepackage{latexsym,euscript}
\usepackage[english,francais]{babel}
\usepackage{times}

\newtheorem{theoremempty}{Theorem}

\newtheorem{maintheorem}[theoremempty]{Main Theorem}

\newtheorem{Theorem}{Theorem}
\newtheorem{theoa}[theoremempty]{Theorem A}
\newtheorem{theob}[theoremempty]{Theorem B}

\newtheorem{Lemma}[Theorem]{Lemma}
\newtheorem{Question}{Question}

 \def\NN{{\mathbb N}}

\def\diff{\operatorname{Diff}}

\def\e{{\varepsilon}}

\ComParit{S0764-4442}
\PIT{FLA}
\PXMA{????}
\Add{?}
\Volume{332}
\Year{2001}
\FirstPage{1}
\LastPage{??}
\AuteurCourant{X. Wang}
\TitreCourant{On the hyperbolicity of C$^1$-generic homoclinic classes} 
\Journal
\Rubrique{Syst\`emes dynamiques}{Dynamical systems}
\PresentePar{First name}{NAME}
\Recu{jour mois ann\'ee}{apr\`es r\'evision}{jour mois ann\'ee}
\begin{document}
\selectlanguage{english}
\title{On the hyperbolicity of C$^1$-generic homoclinic classes
}
\author{%
Xiaodong Wang~$^{\text{a,b}}$
}
\address{%
\begin{itemize}\labelsep=2mm\leftskip=-5mm
\item[$^{\text{a}}$]
School of Mathematical Sciences, Peking University, Beijing, 100871, P.R. China\\
\item[$^{\text{b}}$]
 Laboratoire de Math\'ematiques d'Orsay,
Universit\'e Paris-Sud 11, Orsay 91405, France\\
E-mail: xdwang1987@gmail.com
\end{itemize}
}
\maketitle
\thispagestyle{empty}
\begin{Abstract}{%
Works of Liao, Ma\~n\'e, Franks, Aoki and Hayashi characterized lack of hyperbolicity
for diffeomorphisms by the existence of weak periodic orbits.
In this note we announce a result which can be seen as a local version of these works:
for C$^1$-generic diffeomorphism, a homoclinic class either is hyperbolic
or contains a sequence of periodic orbits that have a Lyapunov exponent arbitrarily close to 0.
}\end{Abstract}
\selectlanguage{french}
\begin{Ftitle}{%
Sur l'hyperbolicit\'e des classes homoclines C$^1$-g\'en\'eriques
}\end{Ftitle}
\begin{Resume}{%
Des travaux de Liao, Ma\~n\'e, Franks, Aoki et Hayashi ont caract\'eris\'e le manque d'hyperbolicit\'e des
diff\'eomorphismes par l'existence d'orbites p\'eriodiques faibles. Dans cette note, nous annon\c{c}ons un
r\'esultat qui peut \^{e}tre consid\'{e}r\'{e} comme une version locale de ces travaux: pour les diff\'{e}omorphismes C$^1$-g\'{e}n\'{e}riques,
une classe homocline ou bien est hyperbolique, ou bien contient une suite d'orbites p\'{e}riodiques qui ont un
exposant de Lyapunov arbitrairement proche de 0.
}\end{Resume}

\par\medskip\centerline{\rule{2cm}{0.2mm}}\medskip
\setcounter{section}{0}
\selectlanguage{english}

\section{Introduction}

It is known for a long time that hyperbolic invariant compact sets have many nice properties, like shadowing properties, the stability properties, the existence of uniform stable and unstable manifolds, etc. So it is important to understand dynamics beyond hyperbolicity and to characterize hyperbolicity. The first works in this direction were devoted to the stability conjecture, which tells that hyperbolic diffeomorphisms are the only ones that are $\Omega$-stable. For surface diffeomorphisms, this conjecture has been solved independently by Liao and Ma\~{n}\'{e} in~\cite{l2} and~\cite{m2}. In their proofs, Liao's selecting lemma and Ma\~{n}\'{e}'s ergodic closing lemma played an important role. For higher dimensions, Ma\~{n}\'{e} solved it in~\cite{m3}. From~\cite{f}, it is easy to see that a C$^1$ $\Omega$-stable diffeomorphism $f$ satisfies the \emph{star condition}: there is a C$^1$-neighborhood $\mathcal{U}$ of $f$, such that any $g\in\mathcal{U}$ has no non-hyperbolic periodic point. Ma\~{n}\'{e} conjectured that if a diffeomorphism satisfied the star condition, then it is hyperbolic, that is to say, its chain recurrent set is hyperbolic. This conjecture was proved by Aoki and Hayashi, see~\cite{ao} and~\cite{h2}. Then one would ask the following question naturally, which is a local version of this conjecture. Recall that a homoclinic class $H(p)$ of a hyperbolic periodic point $p$ is the closure of the union of hyperbolic periodic orbits that are homoclinically related to $orb(p)$.
\medskip

\begin{Question}\emph{(Problem 1.8 in~\cite{csy})}\label{main problem}
For C$^1$-generic $f\in\diff^1(M)$, if a homoclinic class $H(p)$ is not hyperbolic, does it contain periodic orbits that have a Lyapunov exponent arbitrarily close to 0?
\end{Question}
\medskip

The works of Liao, Ma\~n\'e imply the existence of weak periodic orbits close to a non-hyperbolic homoclinic class for C$^1$-generic diffeomorphisms. The difficulty of Question~\ref{main problem} is to link the weak periodic orbits to the homoclinic class. More precisely, if $H(p)$ is not hyperbolic, we can get weak periodic orbits arbitrarily close to it by perturbation with the classical arguments, but we do not know whether they are contained in the homoclinic class. In this paper, we can prove that, generically, they are in fact contained in the homoclinic class.
\medskip

\section{Precise statements}

Let $M$ be a compact connected smooth Riemannian manifold without boundary of dimension $d$. Denote by $\diff^1(M)$ the space of C$^1$-diffeomorphisms from $M$ to $M$. For a diffeomorphism $f\in\diff^1(M)$, for any number $\e>0$, we call a sequence of points $\{x_i\}_{i=a}^b$ an \emph{$\e$-pseudo orbit of $f$}, if $d(f(x_i),x_{i+1})<\e$ for any $i=a,a+1,\cdots,b-1$, where $-\infty\leq a<b\leq \infty$. We say $y$ is \emph{chain attainable} form $x$, denoted by $x\dashv y$, if for any $\e>0$, there is an $\e$-pseudo orbit $\{x=x_0,x_1,\cdots,x_n=y\}$ of $f$. The \emph{chain recurrent set} of a diffeomorphism $f$, denoted by $R(f)$, is the union of points that are chain attainable from itself. We call two points $x,y$ are \emph{chain related}, denoted by $x\sim y$, if $x\dashv y$ and $y\dashv x$. The relation $\sim$ is an equivalent relation on $R(f)$, and every equivalent class of $\sim$ is called a \emph{chain recurrence class}. For a point $x\in R(f)$, denote by $\mathcal{C}(x,f)$ the chain recurrence class that contains $x$.
\medskip

Let $\Lambda$ be an invariant compact set. We say $\Lambda$ is \emph{hyperbolic}, if there is a continuous splitting $T_{\Lambda}M=E^s\oplus E^u$, such that $E^s$ is \emph{contracted} and $E^u$ is \emph{expanded}, that is to say, there are two constants $C>0$ and $\lambda\in(0,1)$, such that, for any $x\in\Lambda$ and any integer $n\in\NN$, we have $\|Df^n|_{E(x)}\|<C\lambda^n$ and $\|Df^{-n}|_{F(x)}\|<C\lambda^n$. If the orbit of a periodic point $p$ is hyperbolic, then we call $p$ a \emph{hyperbolic periodic point}, and the dimension of $E^s$ is called the \emph{index} of $p$, denoted by $ind(p)$.

The set $\Lambda$ is said to have an \emph{dominated splitting}, if there are a continuous splitting $T_{\Lambda}M=E\oplus F$, an integer $m\in\NN$ and a constant $\lambda\in(0,1)$, such that $\|Df^m|_{E(x)}\|\cdot \|Df^{-m}|_{F(f^mx)}\|<\lambda$ for all $x\in\Lambda$. Sometimes we call a dominated splitting associated with the two numbers $m$ and $\lambda$ an \emph{$(m,\lambda)$-dominated splitting}.

If $\mu$ is an ergodic measure for $f$, then there are $d$ numbers $\chi_1\leq\chi_2\leq\cdots\leq\chi_d$, such that for $\mu$ almost point $x\in M$, any non-zero vector $v\in T_xM$, one has $\lim_{n\rightarrow+\infty}\frac{1}{n}\log\|Df^n(v)\|=\chi_i$ for some $i=1,2,\cdots,d$. These numbers are called the \emph{Lyapunov exponents} of the measure $\mu$. Particularly, we call the Lyapunov exponents of the Dirac measure of a periodic orbit the Lyapunov exponents of the periodic orbit. Hence a periodic point $p$ is hyperbolic if and only if all the Lyapunov exponents of $orb(p)$ are non-zero.
\medskip

For any point $x\in M$, any number $\delta>0$, we define the \emph{local stable manifold} and \emph{local unstable manifold} of $x$ of size $\delta$ respectively as follows:\\
$W^s_{\delta}(x)=\{y: d(f^n(x),f^n(y))\leq\delta,\forall n\geq 0; \text{ and } \lim_{n\rightarrow +\infty} d(f^n(x),f^n(y))=0\}$;\\
$W^u_{\delta}(x)=\{y: d(f^{-n}(x),f^{-n}(y))\leq\delta,\forall n\geq 0; \text{ and } \lim_{n\rightarrow +\infty} d(f^{-n}(x),f^{-n}(y))=0\}$.\\
and the \emph{stable manifold} and \emph{unstable manifold} of $x$ respectively as follows:\\
$W^s(x)=\{y:\lim_{n\rightarrow +\infty} d(f^n(x),f^n(y))=0\}$;\\
$W^u(x)=\{y:\lim_{n\rightarrow +\infty} d(f^{-n}(x),f^{-n}(y))=0\}$.\\
By~\cite{hps}, for a hyperbolic invariant compact set $\Lambda$ of $f$, there is a number $\delta>0$, such that for any $x\in\Lambda$, the local stable manifold $W^s_{\delta}(x)$ of $x$ is an embedding disk with dimension $dim(E^s)$ and is tangent to $E^s$ at $x$, where $T_{\Lambda}M=E^s\oplus E^u$ is the hyperbolic splitting. Moreover, the stable manifold $W^s(x)$ of $x$ is an immersing submanifold of $M$. Symmetrically, we have similar statements for $W^u_{\delta}(x)$ and $W^u(x)$.
\medskip

Two hyperbolic periodic points $p$ and $q$ of $f$ are called \emph{homoclinic related}, if their stable and unstable manifolds respectively intersect transversely, that is to say, $W^u(orb(p))\pitchfork W^s(orb(q))\neq \emptyset$ and $W^s(orb(p))\pitchfork W^u(orb(q))\neq\emptyset$. For a hyperbolic periodic point $p$, the closure of the set of periodic points that are homoclinically related to $p$ is called the \emph{homoclinic class} of $p$, denoted by $H(p)$. Also, it is well known that $H(p)$ is the closure of all transverse intersections of its stable and unstable manifolds, that is to say, $H(p)=\overline{W^u(orb(p))\pitchfork W^s(orb(p))}$.
\medskip

For an invariant compact set $\Lambda$ of $f$, a $Df$-invariant sub-bundle $E\subset T_{\Lambda}M$, an integer $m\in\NN$, and any number $\lambda\in(0,1)$, we call $x\in\Lambda$ an \emph{$(m,\lambda)$-$E$-Pliss point}, if $\prod_{i=0}^{n-1} \|Df^{im}|_{E(f^{im}(x))}\|\leq {\lambda}^n$, for any integer $n>0$. If $\Lambda$ does not contain any $(m,\lambda)$-$E$-Pliss point, we call $\Lambda$ an \emph{$(m,\lambda)$-$E$-weak set}. We call two $(m,\lambda)$-$E$-Pliss points $(f^{n}(x),f^{l}(x))$ on a single orbit \emph{consecutive} $(m,\lambda)$-$E$-Pliss points, if $n<l$ and for all $n<k<l$, $f^k(x)$ is not a $(m,\lambda)$-$E$-Pliss point. And if there is a dominated splitting $T_{\Lambda}M=E\oplus F$ on $\Lambda$, we call $x\in \Lambda$ an \emph{$(m,\lambda)$-bi-Pliss point}, if it is an $(m,\lambda)$-$E$-Pliss point for $f$ and an $(m,\lambda)$-$F$-Pliss point for $f^{-1}$. If $m=1$, we will just write $\lambda$-$E$-Pliss point or $\lambda$-$E$-weak set.
\medskip

A subset $R$ of a topological space $X$ is called a \emph{residual} set, if $R$ contains a dense $G_{\delta}$ set of $X$. We say a property is a \emph{generic} property of $X$, if there is a residual set $R\subset X$, such that each element contained in $R$ satisfies the property.
\medskip

We now announce an answer to Question~\ref{main problem}.
\medskip

\begin{maintheorem}
For $C^1$-generic $f\in\diff^1(M)$, a homoclinic class $H(p)$ either is hyperbolic, or contains periodic orbits with arbitrarily long periods that are homoclinically related to $orb(p)$ and have a Lyapunov exponent arbitrarily close to 0.
\end{maintheorem}
\medskip

From~\cite{gy} and Lemma II.3 of~\cite{m2}, we have the fact that for generic $f\in\diff^1(M)$, if all Lyapunov exponents of periodic orbits
that are homoclinically related to $orb(p)$ are uniformly away from 0, then $H(p)$ has a dominated splitting $T_{H(p)}M=E\oplus F$, with $dim E= ind(p)$. Our main theorem is thus a consequence of the following theorem.
\medskip

\begin{theoa}
For $C^1$-generic $f\in \diff^1(M)$, assume that $p$ is a hyperbolic periodic point of $f$. If the homoclinic class $H(p)$ has a dominated splitting $T_{H(p)}M=E\oplus F$, with $dim E\leq ind(p)$, such that the bundle $E$ is not contracted, then there are periodic orbits in $H(p)$ with index $dim(E)$ and with arbitrarily long periods whose maximal Lyapunov exponent along $E$ is arbitrarily close to 0.
\end{theoa}
\medskip

By a standard argument, we can control the norm of product by controlling the product of norm with perturbations. Thus to prove the main theorem, we only have to prove the following.
\medskip

\begin{theob}
For $C^1$-generic $f\in \diff^1(M)$, assume that $p$ is a hyperbolic periodic point of $f$ and that the homoclinic class $H(p)$ has a dominated splitting $T_{H(p)}M=E\oplus F$, with $dim E\leq ind(p)$, such that the bundle $E$ is not contracted. Then there are a constant $\lambda_0\in (0,1)$, and an integer $m_0\in \mathbb{N}$, satisfying: for any $m\in \mathbb{N}$ with $m\geq m_0$, any constants $\lambda_1,\lambda_2\in (\lambda_0,1)$ with $\lambda_1<\lambda_2$, there is a sequence of periodic orbits $\mathcal{O}_k=orb(q_k)$ with period $\tau(q_k)$ contained in $H(P)$, such that
   \begin{displaymath}
        {\lambda_1}^{\tau(q_k)}< \prod_{0\leq i<\tau(q_k)/m} \|Df^m|_{E(f^{im}(q_k))}\|< {\lambda_2}^{\tau(q_k)}.
   \end{displaymath}

\end{theob}

In the next sections, we give a sketch of the proof of the above theorem.

\section{Existence of a bi-Pliss point accumulating backward to an $E$-weak set}

We assume $f$ is a C$^1$-generic diffeomorphism in $\diff^1(M)$ and $H(p)$ is a homoclinic class of $f$ that satisfies the hypothesis of Theorem B. We can choose two numbers $\lambda_0\in (0,1)$ and $m_0\in\mathbb{N}$, such that, for any $m\geq m_0$, the splitting $E\oplus F$ is $(m,\lambda_0)$-dominated, and for the hyperbolic periodic orbit $orb(p)$,
   \begin{displaymath}
     \prod_{0\leq i<\tau(p)/m_0}\|Df^{m_0}|_{E(f^{im_0}(p))}\|<\lambda_0^{\tau(p)/m_0},
   \end{displaymath}
where $\tau(p)$ is the period of $orb(p)$. In the following, we fix $m\geq m_0$. In order to simplify the notations, we will assume that $m=1$ and that $p$ is a fixed point of $f$, but the general case is identical.
\medskip

Since the bundle $E$ is not contracted, there is a point $b\in H(p)$, such that, for any $n\geq 1$, we have $\prod_{i=0}^{n-1} \|Df|_{E(f^{i}(b))}\|\geq 1$. For any number $\lambda\in(\lambda_0,1)$, by Liao's selecting lemma (see ~\cite{l2,w-selecting}), there is a $\lambda$-$E$-weak set contained in $H(p)$ (otherwise, we will get periodic orbits that satisfy the conclusions of Theorem B). Now we fix three numbers $\lambda_1<\lambda_2<\lambda_3\in (\lambda_0,1)$, take the closure of the union of all the $\lambda_2$-$E$-weak sets and denote it by $\hat{K}$. Then there are two cases: either $\hat{K}$ is a $\lambda_2$-$E$-weak set or not. With the arguments related to the Pliss lemma~\cite{p} and the selecting lemma~\cite{l2,w-selecting}, we can get the following lemma under the hypothesis of Theorem B.
\medskip

\begin{Lemma}\label{weak set and pliss point}
There are a $\lambda_2$-$E$-weak set $K\subset H(p)$, a $\lambda_3$-bi-Pliss point $x\in H(p)\setminus K$ satisfying: $\alpha(x)= K$.
\end{Lemma}

\section{The perturbation to make $W^u(p)$ accumulate to the weak set $K$}

Since the $\lambda_2$-$E$-weak set $K$ is contained in $H(p)$, and $W^u(p)$ is dense in $H(p)$, with the technics in the proof of Proposition 10 in~\cite{c2}, we can prove that there is a point on $W^u(p)$ that accumulates the weak set $K$ for a diffeomorphism $g_1$ that is C$^1$ close to $f$ and coincides with $f$ on $K$, $orb(p)$ and the backward orbit of $x$. Moreover, the key point here is that, by the generic assumption of $f$, we can assure that there is a $\lambda_3$-$E$-Pliss point close to $x$ that is contained in the chain recurrence class of $p$ for $g_1$, hence $K$ is still chain related with $p$ for $g_1$. More precisely, we have the following lemma.
\medskip

\begin{Lemma}\label{first perturbation}
There is a residual set $\mathcal{R}\subset \diff^1(M)$, such that, if $f\in\mathcal{R}$ and satisfy the hypothesis of Theorem B, then for any neighborhood $\mathcal{U}$ of $f$ in $\diff^1(M)$, there are a diffeomorphism $g_1\in\mathcal{U}$ and a point $y\in M$, such that,
\begin{enumerate}
\item $y\in W^u(p,f)$ and $g_1$ coincides with $f$ on the set $K\cup orb(p)\cup orb^-(y)$, hence $y\in W^u(p,g_1)$,
\item $\omega(y,g_1)\subset K\subset \mathcal{C}(p,g_1)$.
\end{enumerate}
\end{Lemma}

%

\section{The perturbations to get a heteroclinic connection between $p$ and $K$ }

Denote $K_0=\omega(y,g_1)\subset K$, then for the diffeomorphism $g_1$, $K_0$ is a $\lambda_2$-$E$-weak set, and the orbit of $y$ connects the hyperbolic fixed point $p$ to $K_0$. By two additional perturbations, we can obtain furthermore an orbit connecting $K_0$ to $p$. First, since $K_0\subset \mathcal{C}(p,g_1)$, with the technics for the connecting of pseudo-orbits in~\cite{bc,c1}, we can connect $K_0$ by a true orbit to any neighborhood of $p$ by a $C^1$ small perturbation. Then, by the hyperbolicity of $p$, we use the uniform connecting lemma (see~\cite{w2,wx}) to ``push'' this orbit inside the stable manifold of $p$. In these two steps, the orbit $orb(y)$ that connects $p$ to $K_0$ is not changed. We point out that the proof here is delicate (in fact the most delicate part of the whole proof) and one has to go back in the arguments of ~\cite{bc,c1}.
\medskip

\begin{Lemma}\label{second perturbation}
For the diffeomorphism $f\in\mathcal{R}$, for any neighborhood $\mathcal{U}$ of $f$ in $\diff^1(M)$, there are a diffeomorphism $g_2\in\mathcal{U}$ and two points $y,y'\in M$, such that,
\begin{enumerate}
\item $y\in W^u(p,g_2)$ and $\omega(y,g_2)\subset K$,
\item $y'\in W^s(p,g_2)$ and $\alpha(y',g_2)\subset \omega(y,g_2)$.
\item $g_2$ coincides with $f$ on the set $\omega(y,g_2)\cup orb(p)$,
\end{enumerate}
\end{Lemma}

\section{Last perturbation to get a weak periodic orbit}

Now we have obtained heteroclinic connections between the hyperbolic fixed point $p$ and a subset $K_0=\omega(y,g_2)$ of the weak set $K$. Then using the connecting lemma, we can get a periodic orbit that spends a given proportion of time close to $orb(p)$ and $K_0$ by C$^1$ small perturbation. More precisely, the periodic orbit that we get spends a long time close to the weak set $K_0$, and spends another long time (which can be controlled) close to $p$, hence the average of the product of the norm along the bundle $E$ of this periodic orbit is larger than $\lambda_1$ (controlled by the norm of points close to $K_0$) and smaller than $\lambda_2$ (modified by the norm of points close to $p$). The key point in the connecting process is that, for the hyperbolic fixed point $p$, and the two points $y$ and $y'$, by the $\lambda$-Lemma, there are a number $l$ and two small neighborhoods $U_y$ and $U_{y'}$ of $y$ and $y'$ respectively, such that, for any $n\geq l$, there is an orbit segment with length $n$ that connects $U_{y'}$ to $U_y$, and moreover, only the two endpoints of the segment is contains in $U_{y'}\cup U_y$, and the other part of the segment is close to the point $p$.
\medskip

\begin{Lemma}\label{choose time}
For the diffeomorphism $f\in\mathcal{R}$, for any neighborhood $\mathcal{U}$ of $f$ in $\diff^1(M)$, for any integer $L>0$, any neighborhood $U_p$ of $p$, there is $g\in\mathcal{U}$, such that, $g=f|_{orb(p)}$ and $g$ has a periodic orbit $\mathcal{O}=orb(q)$ with period $\tau>L$ such that $q\in U_p$ and
   \begin{displaymath}
     {\lambda_1}^{\tau}\leq \prod_{0\leq i\leq \tau-1} \|Dg|_{E(g^{i}(q))}\|\leq {\lambda_2}^{\tau}.
   \end{displaymath}
\end{Lemma}
\medskip

Finally, by a standard Baire argument (see for example~\cite{gw}), for the C$^1$-generic diffeomorphism $f$, there is a sequence of periodic orbits that are homoclinically related with each other and accumulates to a subset of $H(p)$, and the product of the norm along the bundle $E$ of these periodic orbits satisfy the inequality in Lemma~\ref{choose time}. Therefore, these periodic orbits are contained in $H(p)$. This finishes the proof of Theorem B.
\medskip


\paragraph{Acknowledgments} The author is grateful to Sylvain Crovisier for many useful suggestions both in solving the problem and in the writing and to Lan Wen for carefully listening to the proof and for the useful comments. The author would also like to thank Shaobo Gan, Dawei Yang and Rafael Potrie for listening to the proof and for the helpful discussions. The author was supported by China Scholarship Council (CSC) 201306010008.

\end{document}